\documentclass{article}
\usepackage{graphicx, amssymb, amsmath}
\usepackage{amsthm}
\usepackage{indentfirst}
\usepackage{mathrsfs}
\usepackage{amsfonts}
\usepackage{color}
\usepackage{graphicx}
\usepackage{algorithm}
\usepackage{algorithmic}
\usepackage{upgreek}
\usepackage{subfigure}

\def\k{\kern .5em}
\def\er{\kern .2em}

\begin{document}
\author{}
\newcommand{\be}{\begin{equation}}
\newcommand{\ee}{\end{equation}}
\newcommand{\ba}{\begin{array}}
\newcommand{\ea}{\end{array}}
\newcommand{\beas}{\begin{eqnarray*}}
\newcommand{\eeas}{\end{eqnarray*}}
\newcommand{\bea}{\begin{eqnarray}}
\newcommand{\eea}{\end{eqnarray}}
\newcommand{\ome}{\Omega}

\newtheorem{theorem}{Theorem}[section]
\newtheorem{lemma}{Lemma}[section]
\newtheorem{remark}{Remark}[section]
\newtheorem{proposition}{Proposition}[section]
\newtheorem{definition}{Definition}[section]
\newtheorem{corollary}{Corollary}[section]

\newtheorem{theo}{Theorem}[section]
\newtheorem{lemm}{Lemma}[section]
\newcommand{\blem}{\begin{lemma}}
\newcommand{\elem}{\end{lemma}}
\newcommand{\bthe}{\begin{theorem}}
\newcommand{\ethe}{\end{theorem}}
\newtheorem{prop}{Proposition}[section]
\newcommand{\bprop}{\begin{proposition}}
\newcommand{\eprop}{\end{proposition}}
\newtheorem{defi}{Definition}[section]
\newtheorem{coro}{Corollary}[section]
\newtheorem{algo}{Algorithm}[section]
\newtheorem{rema}{Remark}[section]
\newtheorem{property}{Property}[section]
\newtheorem{assu}{Assumption}[section]
\newtheorem{exam}{Example}[section]

\renewcommand{\theequation}{\arabic{section}.\arabic{equation}}
\renewcommand{\thetheorem}{\arabic{section}.\arabic{theorem}}
\renewcommand{\thelemma}{\arabic{section}.\arabic{lemma}}
\renewcommand{\theproposition}{\arabic{section}.\arabic{proposition}}
\renewcommand{\thedefinition}{\arabic{section}.\arabic{definition}}
\renewcommand{\thecorollary}{\arabic{section}.\arabic{corollary}}
\renewcommand{\thealgorithm}{\arabic{section}.\arabic{algorithm}}
\newcommand{\lan}{\langle}
\newcommand{\curl}{{\bf curl \;}}
\newcommand{\rot}{{\rm curl}}
\newcommand{\grad}{{\bf grad \;}}
\newcommand{\dvg}{{\rm div \,}}
\newcommand{\ran}{\rangle}
\newcommand{\bR}{\mbox{\bf R}}
\newcommand{\bRn}{{\bf R}^3}
\newcommand{\Coinf}{C_0^{\infty}}
\newcommand{\disp}{\displaystyle}
\newcommand{\ra}{\rightarrow}
\newcommand{\Ra}{\Rightarrow}
\newcommand{\ud}{u_{\delta}}
\newcommand{\Ed}{E_{\delta}}
\newcommand{\Hd}{H_{\delta}}
\newcommand\varep{\varepsilon}
\title{ A Two-grid Method for Linearizing and Symmetrizing the Steady-state Poisson-Nernst-Planck Equations
%\thanks{This work was partially supported
%by the National Science Foundation of China
% under grant 11001062 and the fund from Education Department of Guangxi Province under grant 201012MS094.}
}

\author{ Xuefang Li$^1$
 \and Ying Yang$^{2,*}$
 %\thanks{Corresponding author. School of  Mathematics and Computing Science, Guangxi Colleges and Universities Key Laboratory of Data Analysis and Computation, Guangxi Key Laboratory of Cryptography and information Security, Guilin University of Electronic Technology, Guilin, Guangxi 541004, China. E-mail: yangying@lsec.cc.ac.cn}
 \and Hang Cheng$^3$
}\footnotetext[1]
{ School of  Mathematics and Computing Science, Guilin University of Electronic Technology, Guilin, Guangxi 541004, China.}
\footnotetext[2]
{$^{,*}$Corresponding author. School of  Mathematics and Computing Science, Guangxi Colleges and Universities Key Laboratory of Data Analysis and Computation, Guangxi Key Laboratory of Cryptography and information Security, Guilin University of Electronic Technology, Guilin, Guangxi 541004, China. E-mail: yangying@lsec.cc.ac.cn}
\footnotetext[3]
{ School of  Mathematics and Computing Science, Guilin University of Electronic Technology, Guilin, Guangxi 541004, China.}
\date{}
\maketitle
%\pagenumbering{arabic}

\begin{abstract}
In this paper, a two-grid method is proposed to linearize and symmetrize the steady-state Poisson-Nernst-Planck equations. The computational system is decoupled to linearize and symmetrize equations by using this method, which can improve the computational efficiency compared with the finite element method. Error estimates are derived for the proposed method. The numerical results show that the two-grid method is effective.
  %In this algorithms, a nonlinear and coupled problem is first solved on a coarse grid, and then a linear  and some linearized symmetrized and decoupled two-grid algorithms are proposed based on the decoupled two-grid algorithms for solving the nonlinear nonsymmetric and coupled steady-state Poisson-Nernst-Planck equations. We implement these algorithms and present the numerical results justifying our theoretical error estimates. The numerical tests also show that our algorithms are much more efficient than the standard finite element method and the decoupled two-grid algorithms.
 % Poisson-Nernst-Planck
% equations are a coupled system of nonlinear partial differential equations consisting of the Nernst-Planck equation and
% the electrostatic Poisson equation with delta distribution sources, which describe the electrodiffusion of ions in a solvated
% biomolecular system. In this paper, some error bounds for a piecewise
% finite element approximation to this problem are derived. Several numerical
% examples including biomolecular problems are shown to support our analysis.
\end{abstract}
\noindent {\bf Key words.} Poisson-Nernst-Planck equations, two-grid method, finite element method, linearization, symmetrization, decoupling

\noindent
{\bf 2000 AMS subject classifications.} 65N30, 92C40.

\section{Introduction}{\label{sec1}
%part 1. The background of PNP.
%part 2. Some problems in the computing of PNP
%Part 3. Some ways for decoupling the coupled equation
%part 4. The idea of two-grid method and decoupling
 Electrodiffusion model has aroused a well-accepted publicity over the past years due to its widely applications in semiconductors \cite{j96,m86,n91}, electrochemical systems \cite{bks09,mcr12} and biological membrane channels \cite{bsdk09,ck05,ehl10}. It is mainly described by the Poisson-Nernst-Planck (PNP) equations, which are a system coupled by  the NP equation and the electrostastic Poisson equation. Here the NP equation is nonlinear and nonsymmetric. Due to the nonlinearity and coupling of the PNP equations, in general, it is almost impossible to find the analytic solutions of the equations.
\par  In the past several decades, there appears many literatures on numerical methods for the PNP equations, including finite difference method, finite volume method and finite element method, etc. Finite difference method was widely used to solve the PNP equations \cite{cc65,cck00}, but the main disadvantage of it may be the poor adaptability of irregular region. Finite volume method, which focuses on avoiding the disadvantage of finite difference method, was then applied to solve the PNP equations in irregular domains, but it is not easy to achieve the high accuracy owing to the difficulty of the designing higher-order control volume \cite{mm09,wsx02}. Finite element method, which is suitable for irregular surface, has shown the efficiency and effectiveness of dealing with the electrodiffusion model \cite{lhmz,lzhb}.
\par  The PNP equations are a kind of coupled system. Generally speaking, it is convenient to solve a coupled system by using a decoupling method than solving it directly. The decoupling method in which the coupled problems can be separated into single subproblems, are interested in applications. The suitable methods are allowed to use flexibly for solving each subproblem separately, and the numerical implementation is easy and efficient. Recently, a decoupling two-grid method for the steady-state PNP equations is proposed in \cite{yaluxie}.
By using that method, the coupled system is decoupled by the coarse grid finite element approximation and
good initial values are provided for solving PNP equations, which can improve the computational efficiency
 and save the computational time. The theoretical results show that if the finite element solution on the coarse
grid approximates that on the fine grid well enough, then the two-grid method can achieve the similar
approximation effect that the conventional finite element method could do. The numerical results in \cite{yaluxie} show
the effectiveness of the two-grid method.
\par Note that the two-grid method, proposed originally by Xu \cite{xu92} in 1992, is an efficient numerical method for nonselfadjoint or indefinite problems. Later, Xu \cite{xu94, xu96} designed and analyzed the two-grid method for the nonlinear elliptic equations, in which the main idea is to use a solution on the coarse space to produce a rough approximation of the solution, and then use it as the initial guess so that the original equation is linearized on the fine grid. This procedure involves a nonlinear solver on the coarse space and a linear solver on the fine space. Now this idea has been applied to solving many problems, such as nonlinear fourth-order reaction-diffusion problem \cite{liudu15}, nonlinear parabolic equations \cite{chenchen13} and nonlinear reaction-diffusion equations \cite{wuwu}. In \cite{yaluxie}, the two-grid method is used to decouple the PNP equations. However, the decoupled system is still a nonlinear system or a linear but asymmetric system, so many excellent algorithms which are suitable for the linear and symmetric system can not be used. In fact, we can apply the idea of two-grid method for the linearization and symmetrization of a single equation to the PNP equations, so a two-grid method for linearizing and symmetrizing the PNP equations is generated.

  \par In this paper, we construct two kinds of two-grid algorithms for linearizing and symmetrizing the PNP equations, which are based on the decoupling two-grid algorithms in \cite{yaluxie}. The error estimates are also presented for the algorithms. These error estimates indicate that our two-grid method can retain the same order of approximation accuracy as both the standard finite element method and the two-grid method in \cite{yaluxie} do, but our algorithms can be plugged in the most efficient and optimized local linear and symmetric solvers which are more efficient.
 %save a large amount of CPU running time than the other two numerical methods.
%n these algorithm, the key idea is to use the coarse space to produce a rough approximation of the solution, and then use it as the known-item so that the original equations are linearized, symmetrized and decoupled on the fine grid.

\par The outline of this paper is as follows. Some preliminaries are presented in the next section. %A description and analysis of the finite element approximation are presented in section 3.%
 Two kinds of two-grid algorithms are proposed and analyzed in section 3. Numerical experiments are demonstrated in section 4.

\setcounter{equation}{0}
\section{Preliminaries}\label{sec2}

In this section we aim to describe some notations and
briefly review the PNP equations and their variational forms.
Besides, we introduce the finite element approximation of the PNP equations and some error bounds of them.
%Some numerical results of the decoupled two-grid algorithms \cite{yaluxie} are also shown in this section.
Let
$ \Omega\subset  R^d~(d=2,3)$ be a bounded Lipschitz domain.
We adopt the standard notations for Sobolev spaces $W^{s,p} (\Omega )$
and their associated norms and seminorms \cite{a75}. For $ p = 2$,
the notations $H^s (\Omega ) = W^{s,2} (\Omega )$ and
$H_0^1 (\Omega ) = \{ v \in H^1 (\Omega ):v|_{\partial \Omega } = 0\} ,
$ where $v|_{\partial \Omega } = 0$ is in the sense of trace,
~$\parallel \cdot \parallel_{s,p,\Omega }  = \parallel \cdot \parallel_{W^{s,p} (\Omega )}$ are used.
Let $( \cdot , \cdot )$ denote the standard $L^2$ inner product.
\par We introduce the following steady state PNP equations (see the time-dependent PNP equations in  \cite{suzhli16})
\bea\label{of1}- \nabla  \cdot (\nabla p^i  + q^i p^i \nabla \phi ) = F_i, ~~for ~~i=1,2,\eea
\bea\label{of2} - \Delta \phi  = \sum\limits_{i = 1}^2 {q^i p^i }  + F_3,\eea\\
where $p^i (x):\Omega  \to R_0^ +$ is the concentration of ions,
 $\phi (x):\Omega  \to R$ is the electrostatic potential.
The index i stands for the different ionic species and $ q^i$ is the charge of species i.
For simplicity, in the following we choose
$q^1  = 1,q^2  =  - 1$, $F_i ~(i = 1,2,3)$ indicate the reaction source terms.
In addition, we consider the homogeneous Dirichlet boundary conditions as follows:
~$p^1  = p^2  = \phi  = 0$,~on~$\partial \Omega.$

\par
%Suppose $D^i(x),p^i(x)\in L^{\infty}(\ome_s)$
%and $\phi(x)\in W^{1,\infty}(\Omega_s)$.
The weak formulation of (\ref{of1})-(\ref{of2}) reads: find
$p^i  \in H_0^1 (\Omega ),~i = 1,2$ and $
\phi  \in H_0^1 (\Omega )$ such that %%
\bea\label{wef1}(\nabla p^i,\nabla v) + (q^ip^i\nabla \phi,\nabla v) = (F_i,v),~~\forall v\in
 H_0^1 (\Omega),~~i=1,2,\eea
\bea\label{wef2}(\nabla \phi ,\nabla w) = \sum\limits_{i = 1}^2 {q^i (p^i ,w )}  + (F_3 ,w),
~~\forall w\in H_0^1(\Omega).\eea
%The week formulations of (\ref{prob1}) and (\ref{bd1}) are that:
%Find $\phi\in H^1(\ome_s)$
%and $p^i\in V=\{v|v\in H^1(\ome_s),v|_{\partial\ome}=p^i_{pulk}\}$
%$(1\leq i\leq n)$ such that %%
 %\bea\label{wf1}a_1(p^i,v)+b_1(p^i,\phi,v)=0,~~\forall v\in
%V_0,\eea \bea\label{wf2}a_2(\phi,w)+b_2(p^i,w)=f(w),\forall
%w\in C_0^{\infty}(\ome),\eea where \bea
%\label{a1b1}a_1(p^i,v)=(D^i\nabla p^i,\nabla
%v),~~~~b_1(p^i,\phi,v)=(D^i\beta q^i p^i\nabla \phi,\nabla v), \eea
%\bea \label{a2b2} a_2(\phi,w)=(\epsilon \nabla \phi,\nabla
%w),~~~~b_2(p^i,w)=(-\lambda\sum_{i=1}^n  q^i p^i,w),~~f(w)=(\rho^f,w) \eea and the space
%$V_0=\{v|v\in H^1(\ome_s),v|_{\partial\ome}=0\}$.

%Define a Galerkin-projection
% $R_h: H_0^1(\ome)\rightarrow S_0^h(\ome)$
%  by
%%
%\bea    a(u-R_hu,v)=0,~~ \forall v\in S_0^h(\ome).\eea
%%
%Define a Galerkin-projection $R_h^*:H_0^1(\ome)\rightarrow
%S_0^h(\ome)$ by \beas a(v,u-R_h^*u)=0,~~ \forall v\in
%S_0^h(\ome).\eeas
%
%From \cite{yazh01}, we have the following error estimation
%  \bea\label{l2pro}\|P_hw-R_h^*w\|_{0,\infty,\ome_m}\leq C h^{\frac 1 2}|\ln {\frac d h}| \|w\|_{2,\ome},\forall w\in H_0^1(\ome)\cap W^{2,2}(\ome). \eea
%

 %%%%%%%%%%%%%%%%%%%%%%%%%%%%%%%%%%%%%%%%%%%%%%%%%%%%%%%%%%%%%%%%%%%%%%%%%%%%%5
% \section{The finite element approximation}\label{sec3}\setcounter{equation}{0}

\par Let $\Gamma _h$ be a quasi-uniform triangulation of $\Omega$ with mesh size  $h > 0$
and define corresponding linear finite element space
\bea\label{kongjian}
 S_h(\Omega)  = \{ v \in H^1 (\Omega ):v|_{\partial\ome}=0 ~~\mbox{and}~~ v|_e \in p^1 (e),
\forall e \in \Gamma_h (\Omega )\},\eea
 ~where $p^1 (e)$ is the set of linear polynomials. The coarse space $S_H(\Omega)$ is defined by replacing
$h$ with $H$. Moreover,
throughout this paper C denotes a positive constant independent of $h$,
but may have different values at different places.
\par Assume there exists a unique solution $(\phi ,p^i )~(i = 1,2)$
satisfying (\ref{wef1})-(\ref{wef2}). The standard finite element
approximation of problem (\ref{wef1})-(\ref{wef2}) is defined as flows: find
$(p_h^1 ,p_h^2 ,\phi _h ) \in [S_h(\Omega) ]^3$ such that,
\bea\label{fem1}(\nabla p_h^i,\nabla v_h) + (q^ip_h^i\nabla \phi_h,\nabla v_h) = (F_i,v_h),
 ~~\forall v_h\in S_h (\Omega) ,~~i=1,2,\eea
\bea\label{fem2}(\nabla \phi_h ,\nabla w_h) = \sum\limits_{i = 1}^2 {q^i (p_h^i ,w_h )}  + (F_3 ,w_h),
~~\forall w_h\in  S_h(\Omega).\eea
\par We assume there exists a unique solution $(\phi_h, p^i_h)$ satisfying (\ref{fem1}) and (\ref{fem2}). Some error bounds for the finite element approximations was presented in \cite{yalu13}.
If $
\phi  \in H^{1 + m} (\Omega _s )$ and $p^i  \in H^{1 + m} (\Omega _s )~(1 \le i \le n,0 < m \le 1),$
we can obtain
\bea\label{bf1}\parallel \phi  - \phi _h \parallel_{1,\Omega _s }  \le C(h^m  +\sum\limits_{i = 1}^n {\parallel}p^i  - p_h^i \parallel_{0,\Omega _s } ),\eea
and
\bea\label{bf2}\parallel p^i  - p_h^i \parallel _{1,\Omega _s }  \le C(h^m  +\sum\limits_{i = 1}^n {\parallel} p^i  - p_h^i \parallel_{0,\Omega _s } ),\eea
where~~$p_h^i  \in L^\infty  (\Omega _s )$ is assumed.

%%%%%%%%%%%%%%%%%%%%%%%%%%%%%%%%%%%%%%%%%%%%%%%%%%%%%%%%%%%%%%%%%%%%%%%%%%%%%%%%%%%%%%%%%%%%%%%
\section{The two-grid algorithms}\label{sec3}\setcounter{equation}{0}
In this section, two kinds of two-grid algorithms are designed to linearize and symmetrize the PNP equations.
Some error estimates are also presented for the algorithms.
\par We first recall the decoupling two-grid algorithms in \cite{yaluxie} which provide a basis for our two-grid algorithms.
\begin{algo}\label{algo1}( the decoupling  two-grid Algorithm I \cite{yaluxie})$     $\\
Step 1. On the coarse grid, solve the coupled system (\ref{fem1}) and (\ref{fem2}) as follows: find $(p_H^1,p_H^2,{\phi _H}) \in {[{S_H(\Omega)}]^3}$, such that
\bea\label{dtgc1}(\nabla p_H^i,\nabla {v_H}) + ({q^i}p_H^i\nabla {\phi _H},\nabla {v_H}) = ({F_i},{v_H}),
~~\forall v_H \in S_H (\Omega) ,~~i=1,2,\eea
\bea\label{dtgc2}(\nabla \phi _H ,\nabla w_H ) = \sum\limits_{i = 1}^2 {q^i (p_H^i ,w_H )}  + (F_3 ,w_H ),
~~\forall w_H \in S_H(\Omega).\eea
Step 2. On the fine grid, solve the decoupled system as follows: find $( p_h^{1,*},p_h^{2,*},{\phi _h^*}) \in {[S_h(\Omega) ]^3}$, such that
\bea\label{dtgf1}(\nabla p_h^{i,*},\nabla {v_h}) + ({q^i}p_h^{i,*}\nabla {\phi _h^*},\nabla {v_h}) = ({F_i},{v_h}),
~~\forall v_h \in S_h (\Omega),~~i=1,2,\eea
\bea\label{dtgf2}(\nabla \phi _h^* ,\nabla w_h ) = \sum\limits_{i = 1}^2 {q^i (p_H^i ,w_h )}  + (F_3 ,w_h ),
~~\forall w_h \in S_h(\Omega).\eea
\end{algo}
\par Algorithm \ref{algo1} is a semi-decoupling two-grid algorithm, since the  solution of density in (\ref{dtgf1}) still depends on the solution of potential in (\ref{dtgf2}). A fully decoupling algorithm is also presented in \cite{yaluxie} as follows:
%In the next algorithm, the system (\ref{dtg2f1}) and (\ref{dtg2f2}) are fully decoupled, so it can be solved in parallel on the fine grid %level. That is the difference between Algorithm \ref{algo1} and Algorithm \ref{algo2}.
\begin{algo}\label{algo2}( the decoupling two-grid Algorithm II \cite{yaluxie})$     $\\
Step 1. On the coarse grid, solve the coupled system (\ref{fem1}) and (\ref{fem2}) as follows: find $(p_H^1,p_H^2,{\phi _H}) \in {[{S_H(\Omega)}]^3}$, such that
\bea\label{dtg2c1}(\nabla p_H^i,\nabla {v_H}) + ({q^i}p_H^i\nabla {\phi _H},\nabla {v_H}) = ({F_i},{v_H}),
~~\forall v_H \in S_H (\Omega) ,~~i=1,2,\eea
\bea\label{dtg2c2}(\nabla \phi _H ,\nabla w_H ) = \sum\limits_{i = 1}^2 {q^i (p_H^i ,w_H )}  + (F_3 ,w_H ),
~~\forall w_H \in S_H(\Omega).\eea
Step 2. On the fine grid, solve the decoupled system as follows: find $( p_h^{1,*},p_h^{2,*},{\phi _h^*}) \in {[S_h(\Omega) ]^3}$, such that
\bea\label{dtg2f1}(\nabla p_h^{i,*},\nabla {v_h}) + ({q^i}p_h^{i,*}\nabla {\phi _H},\nabla {v_h}) = ({F_i},{v_h}),
~~\forall v_h \in S_h (\Omega),~~i=1,2,\eea
\bea\label{dtg2f2}(\nabla \phi _h^* ,\nabla w_h ) = \sum\limits_{i = 1}^2 {q^i (p_H^i ,w_h )}  + (F_3 ,w_h ),
~~\forall w_h \in S_h(\Omega),\eea
\end{algo}

\par %In the above decoupled two-grid algorithms, the NP equation is still nonlinear so that the compute speed is too low. In order to improve %the speed of calculation of the algorithms, we now propose two linearized and decoupled two-grid algorithms based on the above decoupled %two-grid algorithms.
Compared with Algorithm \ref{algo1}, the system on the fine grid is fully decoupled in this algorithm, so it can be solved in parallel. Note that we have to solve a kind of Nernst-Planck-like equation on the fine grid which is nonlinear asymmetric in Algorithm \ref{algo1} and linear asymmetric in Algorithm \ref{algo2} (see (\ref{dtgf1}) and (\ref{dtg2f1})). In fact, this kind of equation can be linearized and symmetrized by using the coarse grid solution, which generates the following Algorithm \ref{algo3} and \ref{algo4}.

\begin{algo}\label{algo3}$     $\\
Step 1. On the coarse grid, solve the following nonlinear and coupled system for $(p_H^1,p_H^2,{\phi _H}) \in {[{S_H(\Omega)}]^3}:$
\bea\label{tgfem2c1}(\nabla p_H^i,\nabla {v_H}) + ({q^i}p_H^i\nabla {\phi _H},\nabla {v_H}) = ({F_i},{v_H}),
~~\forall v_H \in S_H (\Omega) ,~~i=1,2,\eea
\bea\label{tgfem2c2}(\nabla \phi _H ,\nabla w_H ) = \sum\limits_{i = 1}^2 {q^i (p_H^i ,w_H )}  + (F_3 ,w_H ),
~~\forall w_H \in S_H(\Omega).\eea
Step 2. On the fine grid, we first solve the Poisson equation for $\phi _h^* \in S_h(\Omega) :$
\bea\label{tgfem2f1}(\nabla \phi _h^* ,\nabla w_h ) = \sum\limits_{i = 1}^2 {q^i (p_H^i ,w_h )}  + (F_3 ,w_h ),
~~\forall w_h \in S_h(\Omega),\eea
then we solve the following linear and symmetric equation on the fine grid for $ (p_h^{1,*},p_h^{2,*}) \in [S_h(\Omega)]^2:$
\bea\label{tgfem2f2}(\nabla p_h^{i,*},\nabla {v_h}) + ({q^i}p_H^i\nabla {\phi _h^*},\nabla {v_h}) = ({F_i},{v_h}),
~~\forall v_h \in S_h (\Omega),~~i=1,2.\eea
\end{algo}
\par
 %In this algorithm, a nonlinear and coupled system is first solved on a coarse grid, then a linear symmetric and decoupled system is solved on a fine grid by using the coarse grid solutions.
 %In Step 1, we adopt the standard finite element method to discretize the equations on the coarse grid space which is similar to Algorithm \ref{algo1}. In Step 2, the system is  semi-decoupled.
 %The advantages of the semi-decoupled system are introduced in \cite{yaluxie}, for example, the semi-decoupled system could %reduce plenty of computational time.
The Nernst-Planck equation is reduced to a linear and symmetric one on the fine discrete grid by using the coarse grid solution of concentration, which is the difference between Algorithm \ref{algo3} and Algorithm \ref{algo1} (see (\ref{tgfem2f2}) and (\ref{dtgf1})). Since the system on the fine grid is linear symmetric and decoupled, the equations
can be solved separately and efficiently by using the existing optimized computing softwares for linear symmetric equations. Next the error analysis for Algorithm \ref{algo3} is presented.

\begin{theo}\label{theo1}Let $(\ p_h^i, \phi_h)$,~ $(\ p_H^i,\phi_H)$ and $(\ p_h^{i,*},\phi_h^*)$ be the solutions of (\ref{fem1})-(\ref{fem2}), (\ref{tgfem2c1})-(\ref{tgfem2c2}) and (\ref{tgfem2f1})-(\ref{tgfem2f2}) respectively,
then there holds
\bea\label{t2eq1}\parallel \nabla (\phi _h  - \phi _h^* )\parallel _{0,\Omega } \le C\sum\limits_{i = 1}^2\parallel p_h^i  - p_H^i \parallel _{0,\Omega }.\eea
\end{theo}
\textit{Proof.} Comparing (\ref{fem2}) with (\ref{tgfem2f1}), we have
$$(\nabla ({\phi _h} - \phi _h^*),\nabla {w_h}) = \sum\limits_{i = 1}^2 {{q^i}} (p_h^i  - p_H^i,{w_h}), ~~\forall{w_h} \in {S_h(\Omega)}.$$
Let ${w_h} = {\phi _h} - \phi _h^*$ and by Poincar\'e inequality, we can obtain
\beas\parallel \nabla (\phi _h  - \phi _h^* )\parallel _{0,\Omega }^2  &= & |\sum\limits_{i = 1}^2 {q^i } (p_h^i  - p_H^i ,\phi _h  - \phi _h^* )|\nonumber\\& \le & C \sum\limits_{i = 1}^2  \parallel p_h^i  - p_H^i \parallel _{0,\Omega } \parallel \phi _h  - \phi _h^* \parallel _{0,\Omega } \nonumber\\& \le & C\sum\limits_{i = 1}^2\parallel p_h^i  - p_H^i \parallel _{0,\Omega } \parallel \nabla (\phi _h  - \phi _h^*) \parallel _{0,\Omega }.\eeas Hence, we get$$\parallel \nabla (\phi _h  - \phi _h^* )\parallel _{0,\Omega } \le C\sum\limits_{i = 1}^2\parallel p_h^i  - p_H^i \parallel _{0,\Omega }.$$ This completes the proof.$\hfill\Box$

\begin{theo}\label{theo2}Let $(\ p_h^i, \phi_h)$,~ $(\ p_H^i,\phi_H)$ and $(\ p_h^{i,*},\phi_h^*)$ be the solutions of (\ref{fem1})-(\ref{fem2}), (\ref{tgfem2c1})-(\ref{tgfem2c2}) and (\ref{tgfem2f1})-(\ref{tgfem2f2}) respectively, and suppose $\phi _h  \in W^{1,\infty } (\Omega )$~and~$p_H^i  \in L^\infty  (\Omega )$,~ then we have
\bea\label{t2eq2}\parallel  \nabla( p_h^i  - p_h^{i,*} )\parallel _{0,\Omega }  \le C\sum\limits_{i = 1}^2\parallel p_h^i  - p_H ^i\parallel _{0,\Omega }.\eea
\end{theo}
\textit{Proof.} Subtracting (\ref{fem1}) from (\ref{tgfem2f2}), we obtain
$$(\nabla (p_h^i  - p_h^{i,*}),\nabla {v_h}) + ({q^i}{p_h^i}\nabla {\phi _h} - {q^i}{p_H^i}\nabla \phi _h^*,\nabla {v_h}) = 0, ~~\forall{v_h} \in {S_h(\Omega)}.$$
Setting ${v_h} = p_h^i  - p_h^{i,*}$, we get
\bea\label{peq1}\parallel \nabla (p_h^i  - p_h^{i,*} )\parallel _{0,\Omega }^2  & = &  |({q^i}{p_h^i}\nabla {\phi _h} - {q^i}{p_H^i}\nabla \phi _h^*,\nabla (p_h^i  - p_h^{i,*})|\nonumber\\& \le & C \parallel{p_h^i}\nabla {\phi _h}-{p_H^i}\nabla \phi _h^* \parallel _{0,\Omega } \parallel \nabla (p_h^i  - p_h^{i,*})\parallel _{0,\Omega }.\eea
If $\phi _h  \in W^{1,\infty } (\Omega )$ and $p_H^i  \in L^\infty  (\Omega ),$ we have
\bea\label{peq2} \parallel{p_h^i}\nabla {\phi _h}-{p_H^i}\nabla \phi _h^* \parallel _{0,\Omega }
& \le &\parallel{p_h^i}\nabla {\phi _h}-{p_H^i}\nabla {\phi _h} \parallel _{0,\Omega }+\parallel{p_H^i}\nabla {\phi _h}-{p_H^i}\nabla {\phi _h^*} \parallel _{0,\Omega }\nonumber\\&\le & \parallel{p_h^i}-{p_H^i}\parallel _{0,\Omega }+\parallel\nabla( \phi _h- \phi _h^*) \parallel _{0,\Omega }\nonumber\\&\le&C
\sum\limits_{i = 1}^2\parallel{p_h^i}-{p_H^i}\parallel _{0,\Omega },\eea where we have used Theorem \ref{theo1}. Inserting (\ref{peq2}) into (\ref{peq1}), we can easily get the result of Theorem \ref{theo2}. This completes the proof.$\hfill\Box$

\begin{remark} \label{rema1}Theorem \ref{theo1} and Theorem \ref{theo2} indicate that the $H^1$ norm errors between the finite element solutions and the solutions in Algorithm \ref{algo3} are governed by the errors between the concentration on the coarse grid  $(p_H^i)$ and the concentration on the fine grid $(p_h^i)$. If we assume $p_h^i$ approximates to $p_H^i$ well enough, for instance, \beas \sum\limits_{i = 1}^2\parallel{p_h^i}-{p_H^i}\parallel _{0,\Omega } = O(H^2 ),\eeas%then the error
%$\sum_{i=1}^n\|p_h^i-p_H^i\|_{0,\ome_s}$ may achieve optimal convergence
%rate, i.e.
then we can conclude that the Algorithm \ref{algo3} can retain the same order of approximation accuracy as the standard finite element method does under the condition $h = H^2$, that is \bea
\parallel \nabla ( \phi_h - \phi_h^* )\parallel _{0, \Omega}=O(H^2)=O(h)\eea and \bea \parallel
\nabla(p_h^i-p_h^{i,*})\parallel _{0,\Omega} =O(H^2)=O(h).\eea
\end{remark}

We will demonstrate the second algorithm as follows.
%%%%%%%%%%%%%%%%%%%%%%%%
%%%%%%%%%%%%%%%%%%%%%%%% ÒÔÏ ˫ÐÐ %% ¼Ðס²¿·ÖÊÇмӽøÈ¥µÄËã·¨ 3.4 ¼°ÆäÎó²î¹À¼Æ½á¹û
\begin{algo}\label{algo4}$     $\\
Step 1. On the coarse grid, solve the following nonlinear and coupled system for~$(p_H^1,p_H^2,{\phi _H}) \in {[{S_H(\Omega)}]^3}:$
\bea\label{tgfem3c1}(\nabla p_H^i,\nabla {v_H}) + ({q^i}p_H^i\nabla {\phi _H},\nabla {v_H}) = ({F_i},{v_H}),
~~\forall v_H \in S_H (\Omega) ,~~i=1,2,\eea
\bea\label{tgfem3c2}(\nabla \phi _H ,\nabla w_H ) = \sum\limits_{i = 1}^2 {q^i (p_H^i ,w_H )}  + (F_3 ,w_H ),
~~\forall w_H \in S_H(\Omega).\eea
Step 2.  On the fine grid, solve the following linear, symmetric and decoupled system for $
(p_h^{1,*} ,p_h^{2,*} ,\phi _h^* ) \in [S_h(\Omega)]^3:$
\bea\label{tgfem3f1}(\nabla p_h^{i,*},\nabla {v_h})= ({F_i},{v_h})-({q^i}p_H^i\nabla {\phi _H},\nabla {v_h}),
~~\forall v_h \in S_h (\Omega) ,~~i=1,2,\eea
\bea\label{tgfem3f2}(\nabla \phi _h^* ,\nabla w_h ) = \sum\limits_{i = 1}^2 {q^i (p_H^i ,w_h )}  + (F_3 ,w_h ),
~~\forall w_h \in S_h(\Omega).\eea
\end{algo}
\par In this algorithm, the system in Step 2 is  a fully-decoupled one which is the difference between Algorithm \ref{algo3} and Algorithm \ref{algo4}. We can treat these two small linear problem (\ref{tgfem3f1}) and (\ref{tgfem3f2}) in parallel which can save a large amount of CPU running time. The difference between Algorithm \ref{algo4} and Algorithm \ref{algo2} is that the system on the fine grid is a symmetric one for the former but an asymmetric one for the latter.

\par We can get the following error estimates for the solutions of the above algorithm.
\begin{theo}\label{theo3}If $(\ p_h^i, \phi_h)$, $(\ p_H^i,\phi_H)$ and $(\ p_h^{i,*},\phi_h^*)$ are the solutions of (\ref{fem1})-(\ref{fem2}), (\ref{tgfem3c1})-(\ref{tgfem3c2}) and (\ref{tgfem3f1})-(\ref{tgfem3f2}) respectively,
 then we have
\bea\label{t3eq1}\parallel\nabla (\phi _h  - \phi _h^* )\parallel_{0,\Omega } \le C\sum\limits_{i = 1}^2\parallel p_h^i  - p_H^i \parallel_{0,\Omega }.\eea
\end{theo}
\par The proof of this theorem is the same as Theorem \ref{theo1}, since the only difference between Algorithm \ref{algo3} and Algorithm \ref{algo4} is (\ref{tgfem3f1}) which is not used in this proof.

\begin{theo}\label{theo4}If $(\ p_h^i, \phi_h)$, $(\ p_H^i,\phi_H)$ and $(\ p_h^{i,*},\phi_h^*)$ are the solutions of  (\ref{fem1})-(\ref{fem2}), (\ref{tgfem3c1})-(\ref{tgfem3c2}) and (\ref{tgfem3f1})-(\ref{tgfem3f2}) respectively, $\phi _h  \in W^{1,\infty } (\Omega )$ and $p_H^i  \in L^\infty  (\Omega )$, then we have
\bea\label{t3eq2}\parallel \nabla (p_h^i  - p_h^{i,*} )\parallel_{0,\Omega }  \le C(\parallel p_h ^i - p_H ^i\parallel_{0,\Omega }+
\parallel\nabla (\phi _h  - \phi _H )\parallel_{0,\Omega }).\eea
\end{theo}
\textit{Proof.} To derive error estimates, we subtract the equation (\ref{fem1}) from (\ref{tgfem3f1}), we obtain
$$(\nabla (p_h^i  - p_h^{i,*}),\nabla {v_h}) + ({q^i}{p_h^i}\nabla {\phi _h} - {q^i}{p_H^i}\nabla \phi _H,\nabla {v_h}) = 0, ~~\forall{v_h} \in {S_h(\Omega)}.$$
Similarly, taking ${v_h} = p_h^i  - p_h^{i,*}$, then
\bea\label{p2eq1}\parallel \nabla (p_h^i  - p_h^{i,*} )\parallel _{0,\Omega }^2  & = &  |({q^i}{p_h^i}\nabla {\phi _h} - {q^i}{p_H^i}\nabla \phi _H,\nabla (p_h^i  - p_h^{i,*}))|\nonumber\\& \le & C \parallel{p_h^i}\nabla {\phi _h}-{p_H^i}\nabla \phi _H \parallel _{0,\Omega } \parallel \nabla (p_h^i  - p_h^{i,*})\parallel _{0,\Omega }.\eea

Here \bea\label{p2eq2} \parallel{p_h^i}\nabla {\phi _h}-{p_H^i}\nabla \phi _H \parallel _{0,\Omega }
& \le &\parallel{p_h^i}\nabla {\phi _h}-{p_H^i}\nabla {\phi _h} \parallel _{0,\Omega }+\parallel{p_H^i}\nabla {\phi _h}-{p_H^i}\nabla {\phi _H} \parallel _{0,\Omega }\nonumber\\&\le & C (\parallel{p_h^i}-{p_H^i}\parallel _{0,\Omega }+\parallel\nabla( \phi _h- \phi _H) \parallel _{0,\Omega }),\eea where the assumptions $\phi _h  \in W^{1,\infty } (\Omega )$ and $p_H^i  \in L^\infty  (\Omega )$ are used.
 Inserting (\ref{p2eq2}) into (\ref{p2eq1}), we can easily complete the proof of Theorem \ref{theo4}.$\hfill\Box$
\begin{remark} Under the assumption of \beas \sum\limits_{i = 1}^2\parallel{p_h^i}-{p_H^i}\parallel _{0,\Omega } = O(H^2 ),\eeas we have
\bea\parallel \nabla ( \phi_h - \phi_h^* )\parallel _{0, \Omega}=O(H^2)\eea from Theorem \ref{theo3}. That means the convergence rate of electrostatic potential in Algorithm \ref{algo4} is consistent with the finite element method if $h=H^2$ is satisfied. Theorem \ref{theo4} suggests that the two-grid solution for the concentration in Algorithm \ref{algo4} has the same accuracy as the finite element solution only when $h=H$, whereas the experiment results in section 4 demonstrate the optimal convergence rate even under the condition of $h \ne H$. This implies the error estimates in Theorem \ref{theo4} may not be the optimal one.
\end{remark}

% ×¢£ºÕâÀàËÆÓëËã·¨ 3.3£¬µ±Óõ½´ÖÍø¸ñ½â {{\phi _H}} ʱ£¬Ë«³ß¶È½âÀíÂÛÉÏ´ï²»µ½ÓëÓÐÏÞԪͬ½×
%%%%%%%%%%%%%%%
%%%%%%%%%%%%%%% ÒÔÉÏ Ë«ÐÐ %% ¼Ðס²¿·ÖÊÇмӽøÈ¥µÄËã·¨ 3.4 ¼°ÆäÎó²î¹À¼Æ½á¹û

%%%%%%%%%%%%%%%%%%%%%%%%%%%%%%%%%%%%%%%%%%%%%%%%%%%%%%%%%%%%%%%%%%%%%%%%%%%%%%
\section{Numerical experiment}\setcounter{equation}{0}
In this section, we perform a numerical test to show the effectiveness of our algorithms and demonstrate the error estimates we have presented. All the programs are debugged on the Fortran Power Station 4.0 compiler and all results are generated by the same microcomputer.
\par For simplicity, the domain $\Omega  = [0,1]^3$ is specified and we suppose that the boundary conditions are homogeneous. Furthermore, the right hand side functions are set such that the exact solution $(\phi ,p^1 ,p^2 )$ of (\ref{of1})-(\ref{of2})~is as follows
\beas\begin{array}{*{20}c}
   \begin{cases}
   \phi  = \sin \pi x\sin \pi y\sin \pi z  \\
    p^1  = \sin 2\pi x\sin 2\pi y\sin 2\pi z  \\
    p^2  = \sin 3\pi x\sin 3\pi y\sin 3\pi z.  \\
  \end{cases}
   \end{array}\eeas
we use piecewise linear finite elements on the tetrahedral mesh to discretize the equations.
\par To implement Algorithm \ref{algo3}, first the Gummel iteration is used on the coarse
grid. Given the initial value $(p^{1,0},p^{2,0})\in [S_H^0 (\Omega)]^2$, for $m\geq 0$ find $(p^{1,m+1},p^{2,m+1},\phi^{m+1})\in [S_H^0(\Omega)]^3$
such that

\bea\label{giter22}(\nabla\phi^{m+1},\nabla w)=(F_3,w)+\sum_{i=1}^2q^i(p^{i,m},w),\forall
w\in S_H^0(\ome).\eea
\bea\label{giter11}(\nabla p^{i,m+1},\nabla v)+(q^ip^{i,m+1}\nabla\phi^{m+1},\nabla v)=(F_i,v),~~\forall v\in
S_H^0(\ome),~~i=1,2.\eea
The stopping criterion for this iteration is $\|\phi^{m+1}-\phi^m\|_0<10^{-5}$. Suppose $p
^i_H, i = 1, 2$ is the final solution of the concentration in the above iteration. On the second step, the linear, symmetric and decoupled system is solved by using the coarse grid solution as follows: Find $\phi _h^* \in S_h(\Omega)$, such that
\bea\label{tgfem2f11}(\nabla \phi _h^* ,\nabla w_h ) = \sum\limits_{i = 1}^2 {q^i (p_H^i ,w_h )}  + (F_3 ,w_h ),
~~\forall w_h \in S_h(\Omega),\eea
and find $ (p_h^{1,*},p_h^{2,*}) \in [S_h(\Omega)]^2:$
\bea\label{tgfem2f22}(\nabla p_h^{i,*},\nabla {v_h}) + ({q^i}p_H^i\nabla {\phi _h^*},\nabla {v_h}) = ({F_i},{v_h}),
~~\forall v_h \in S_h (\Omega),~~i=1,2.\eea
The errors between the exact solutions and two-grid solutions of Algorithm \ref{algo3} are shown in Table 1. The errors in $H^1$ norm approximate the first-order reduction as $h=H^2$ becomes smaller, which coincides with the theoretical results.
\begin{table}[htbp]
\centering\label{table 1}
\begin{tabular}{|*{6}{c|}}
\hline H&h & $\|\phi_h^*-\phi\|_1$ & $\|p_h^{1,*}-p^1\|_1$&$\|p_h^{2,*}-p^2\|_1$& CPU(S)\\
\hline
 1/2& 1/4& 9.15E-01  &3.03E+00&5.40E+00& 0.06\\ \hline
 1/4& 1/16  &2.44E-01&9.79E-01&2.12E+00&1.89\\ \hline
 1/8 & 1/64& 6.22E-02&2.58E-01&5.69E-01& 503.75\\ \hline
\end{tabular}
\caption{The $H^1$ norm errors between the exact solutions and the two-grid solutions of Algorithm \ref{algo3}
%(Example \ref{ch4ex1}) \label{ch4ex1:t1}
}
\end{table}

\par The two-grid solutions of Algorithm \ref{algo4} are obtained by using the similar computational procedure
in (\ref{giter22})-(\ref{tgfem2f11}), but (\ref{tgfem2f22}) is replaced with the following equation:
\beas\label{tgfem3f1}(\nabla p_h^{i,*},\nabla {v_h})= ({F_i},{v_h})-({q^i}p_H^i\nabla {\phi _H},\nabla {v_h}),
~~\forall v_h \in S_h (\Omega) ,~~i=1,2.\eeas
The errors between the exact solutions and the two-grid solutions of Algorithm \ref{algo4} are shown in Table 2.  The error for the potential $\phi_h^*$ in $H^1$ norm approximates the first-order reduction as $h=H^2$ becomes smaller, which coincides with the theoretical result of Theorem \ref{theo3}. However,
the errors of solutions $p_h^{1,*}$ and $p_h^{2,*}$ can not approximate the first-order reduction respectively under the condition of $h = H^2$, especially when $H=1/8$ and $h=1/64$. We consider this phenomenon is due to a low order approximation for  $\nabla \phi _H$ to $\nabla \phi_h$ (see Theorem \ref{theo4}). This problem can be remedied by using a smaller size of coarse grid (see the result of last row in Table 2).
%but the error of the solutions $p_h^{1,*},~p_h^{2,*}$ can not achieve the second-order reduction, when $h=H^2$ (Comparing Table 4 with Table 1). Such an order reduction may be caused by the low order approximation for potential $\phi_H$ on coarse grid to $\phi_h$ on fine grid (See Theorem \ref{theo5}). This problem can be remedied by using a smaller size of coarse grid (see the result of last row in Table 4). That means two-grid Algorithm 3.3 can also achieve the same order of accuracy as the finite element method could do, if we choose a suitable coarse mesh size $H$.
\begin{table}[htbp]
\centering\label{table 1}
\begin{tabular}{|*{6}{c|}}
\hline H&h & $\|\phi_h^*-\phi\|_1$ & $\|p_h^{1,*}-p^1\|_1$&$\|p_h^{2,*}-p^2\|_1$& CPU(S)\\
\hline
 1/2& 1/4& 9.15E-01  &3.03E+00&5.40E+00& 0.04\\ \hline
 1/4& 1/16  &2.44E-01&9.89E-01&2.12E+00&1.81\\ \hline
 1/8 & 1/64& 6.22E-02&2.91E-01&5.80E-01& 500.06\\ \hline
 1/32 &1/64& 6.09E-02&2.46E-01&5.48E-01& 629.04\\ \hline
\end{tabular}
\caption{The $H^1$ norm errors between the exact solutions and the two-grid solutions of Algorithm \ref{algo4}
%(Example \ref{ch4ex1}) \label{ch4ex1:t1}
}
\end{table}

Next, we will compare Algorithm \ref{algo3} and \ref{algo4} with the finite element method. The $L^2$ norm and $H^1$ norm errors between the exact solutions and the finite element solutions are shown in Table 3 (cf. \cite{yaluxie}). Compared Table 1 with Table 3, we find that the solutions of Algorithm \ref{algo3} remain the same order of accuracy as
the finite element solutions but require much less computational time. For the full decoupled Algorithm \ref{algo4},  if we choose a suitable coarse mesh size $H$, we can achieve the similar effect as the finite element method could do by comparing Table 2 with Table 3. In addition, the CPU running time used by Algorithm \ref{algo4} is much less than that of the finite element method, which indicates the effectiveness of Algorithm \ref{algo4}. Note that the CPU time can be saved much more if Algorithm \ref{algo4} is implemented on the parallel computers.

\begin{table}[htbp]
\centering
\begin{tabular}{|*{7}{c|}}
\hline h &$\|p_h^{1}-p^1\|_0$ & $\|p_h^{1}-p^2\|_0$ & $\|\phi_h-\phi\|_1$ & $\|p_h^{1}-p^1\|_1$&$\|p_h^{2}-p^2\|_1$& CPU(S)\\
\hline
 1/4& 2.41E-01&3.26E-01&9.14E-01&3.03E+00&5.39E+00& 1.5\\ \hline
 1/8& 8.99E-02&1.72E-01&4.80E-01&1.82E+00&3.75E+00& $-$\\ \hline
 1/16  & 2.53E-02&5.59E-02&2.43E-01&9.57E-01&2.10E+00&7.56\\ \hline
  1/32  &6.51E-03&1.50E-02&1.22E-01&4.85E-01&1.09E+00&$-$\\ \hline
 1/64& 1.64E-03&3.83E-03&6.09E-02&2.44E-01&5.47E-01& 2991.73\\ \hline

\end{tabular}
\caption{ The $L^2$ norm and $H^1$ norm errors between the exact solutions and the finite element solutions (cf. \cite{yaluxie}).
%(Example \ref{ch4ex1}) \label{ch4ex1:t1}
}
\end{table}

\par Finally, in order to compare the two-grid Algorithms \ref{algo3} and \ref{algo4} with the two-grid Algorithm \ref{algo1} and \ref{algo2} presented in \cite{yaluxie}, Table 4 and Table 5 list the $H^1$ norm errors between the exact solutions and the solutions of Algorithm 3.1 and 3.2, respectively, where the numerical results are generated by the same environment as Algorithms \ref{algo3} and \ref{algo4}. Both Table 2 and Table 4 show the results of the semi-decoupling algorithms. Table 3 and Table 5 display the results of the full-decoupling algorithms. A comparison of Table 2 with Table 4 shows that there is no big difference in the solutions of Algorithm 3.3 and Algorithm 3.1. However, Algorithm 3.3 is more efficient than Algorithm 3.1 on the CPU running time. Compared Table 3 with Table 5, we can get a similar conclusion for Algorithm 3.4 and 3.2 as the comparison of  Algorithm 3.3 with Algorithm 3.1. The improvement of efficiency of Algorithm 3.3 and 3.4 is due to the use of linear and symmetric solver for the system while a nonlinear or an asymmetric solver is used in Algorithm 3.1 and 3.2.

\begin{table}[htbp]
\centering\label{table1}
\begin{tabular}{|*{6}{c|}}
\hline H&h=$H^2$ & $\|\phi_h^*-\phi\|_1$ & $\|p_h^{1,*}-p^1\|_1$&$\|p_h^{2,*}-p^2\|_1$& CPU(S)\\
\hline
 1/2& 1/4& 9.15E-01&3.03E+00&5.39E+00& 1.2\\ \hline
 1/4& 1/16  &2.44E-01&9.57E-01&2.10E+00&2.2\\ \hline
 1/8 & 1/64& 6.22E-02&2.44E-01&5.47E-01& 830\\ \hline
\end{tabular}
\caption{The $H^1$ norm errors between the exact solutions and the two-grid solutions of Algorithm \ref{algo1} (cf. \cite{yaluxie})
%(Example \ref{ch4ex1}) \label{ch4ex1:t1}
}
\end{table}
\begin{table}[htbp]
\centering\label{table1}
\begin{tabular}{|*{6}{c|}}
\hline H&h & $\|\phi_h^*-\phi\|_1$ & $\|p_h^{1,*}-p^1\|_1$&$\|p_h^{2,*}-p^2\|_1$& CPU(S)\\
\hline
 1/2& 1/4& 9.15E-01  &3.03E+00&5.39E+00& 1.2\\ \hline
 1/4& 1/16  &2.44E-01&9.89E-01&2.10E+00&2.2\\ \hline
 1/8 & 1/64& 6.22E-02&2.92E-01&5.70E-01& 830\\ \hline
 1/32 &1/64& 6.09E-02&2.46E-01&5.48E-01&1189\\ \hline
\end{tabular}
\caption{The $H^1$ norm errors between the exact solutions and the two-grid solutions of Algorithm \ref{algo2} (cf. \cite{yaluxie})
%(Example \ref{ch4ex1}) \label{ch4ex1:t1}
}
\end{table}

\section{Conclusion}
%%------------------------------
In this paper, we construct and analyze a two-grid method for linearizing and symmetrizing the steady-state PNP equations. This method can reduce a nonlinear and coupled system into a linear, symmetric and decoupled system by solving a nonlinear and coupled system on a much smaller space. The numerical results verify the theoretical results and demonstrate the efficiency of the proposed method. Compared to the finite element method and the two-grid method in \cite{yaluxie}, the new two-grid method can not only save a large amount of CPU running time, but also retain the same order of approximation accuracy as they do. In  the future, it is promising to extend this approach to the time-dependent PNP equations and other kinds of modified PNP equations.\par
{\sc Acknowledgement.}
Y. Yang was supported by the China NSF (NSFC 11561016, NSFC 11561015)
and the fund from Education Department of Guangxi Province under
grant (2014GXNSFAA118004).

%%%%%%%%%%%%%%%%%%%%%%%%%%%%%%%%%%%%%%%%%%%%%%%%%%%%%%%%%%%%%%%%%%%%%%


\begin{thebibliography}{999}
\bibliographystyle{plain}
%\bibitem{ad75}R. A. Adams, Sobolev Spaces, Academic Press, New York, 1975.

%\bibitem{bach97} V. Barcilon, D. P. Chen, R. S. Eisenberg, and J. W. Jerome,
%Qualitative properties of steady-state Poisson-Nernst-Planck
%systems: perturbation and simulation study, SIAM J. Appl. Math.,
%{3}, 1997, 631-648.

%\bibitem{bosa09}D. S. Bolintineanu, A. Sayyed-Ahmad, H. T. Davis, and Y. N. Kaznessis, Poisson-NernstPlanck
%models of nonequilibrium ion electrodiffusion through a protegrin transmembrane pore,
%PLOS Comp. Bio., {5}, 2009, e1000277.
%\bibitem{brsc}S. C. Brenner and L. R. Scott, The Mathematical Theory of Finite Element Methods,
%Springer-Verlag, New York, 1994.

%\bibitem{camu09} M. Cai, M. Mu and J. Xu, Numerical solution to a mixed Navier-Stokes/Darcy model by the two-grid approach, SIAM J. Numer. Anal., {47}, 3325-3338.

%\bibitem{caco00}A. E. Cardenas, R. D. Coalson and M. G. Kurnikova,
%Three-dimensional Poisson-Nernst-Planck theory studies: influence of
%membrane electrostatics on gramicidin a channel conductance,
%Biophys. J., {79} (1), 2000, 80-93.

%\bibitem{coco65}  H. Cohen and J. W. Cooley, The numerical solution of the time-dependent Nernst-Planck
%equations, Biophys. J., {5}, 1965, 145-162.

%\bibitem{eich93}R. Eisenberg and D. P. Chen, Poisson-Nernst-Planck(PNP)
%theory of an open ionic channel, Biophys. J., {64} (2), 1993,
%A22-A22.
%\bibitem{eihy10} B. Eisenberg, Y. Hyon and C. Liu, Energy variational analysis of ions in water and channels: Field theory
%for primitive models of complex ionic fluids, J Chem Phys., {133}, 2010, 104104.
%\bibitem{glpa97} R. Glowinski, T. Pan, and J. Periaux, A Lagrange multiplier/fictitious domain method for
%the numerical simulation of incompressible viscous flow around moving grid bodies: I. Case
%where the rigid body motions are known a priori, C. R. Acad. Sci. Paris S¡äer. I Math., {324},
%1997, 361-369.

%\bibitem{imro02} W. Im and B. Roux, Ion permeation and selectivity of OmpF porin a theoretical study based
%on molecular dynamics, brownian dynamics, and continuum electrodiffusion theory, J. Mol.
%Biol., {322}, 2002, 851-869.

%\bibitem{je85}J. W. Jerome, Consistency of semiconductor modeling: an
%existence/stability analysis for the stationary van Boosbroeck
%system, SIAM J. Appl. Math., {45}, 1985, 565-590.
%\bibitem{jeke91}J. W. Jerome and T. Kerkhoven, A finite element approximation theory for the drift-diffusion
%semiconductor model, SIAM J. Numer. Anal., {28}, 1991, 403-422.
%\bibitem{jish06} J. Jin, S. Shu and J. Xu, a two-grid discretization method for decoupling systems of partial differential equations, Math. Comp., {75}, 2006, 1617-1626.
%\bibitem{kuco99}M. G. Kurnikova, R. D. Coalson, P. Graf and A. Nitzan, A
%lattice relaxation algorithm for three-dimensional
%Poisson-Nernst-Planck theory with application to ion transpoort
%through the gramicidin a channel, Biophys. J., {76} (2), 1999,
%642-656.

%\bibitem{luho10}B. Z. Lu, M. J. Holst, J. A. McCammo and Y. C. Zhou, Poisson-Nernst-Planck equations for simulating biomolecular
%diffusion-reaction processes I: finite element solutions, J. Comput. Phys., {229}, 2010, 6979-6994.

%\bibitem{luzh07}B. Z. lu, Y. C. Zhou, G. A. Huber, S. D. Bond,
%M. J. Holst and J. A. McCammon, Electrodiffusion: a continuum modeling
%%framework for biomolecular systems with realistic spatiotemporal
%resolution, J. Chem. Phys., {127}, 2007, 135102.
%\bibitem{luzh11}B. Z. lu, Y. C. Zhou, Poisson-Nernst-Planck Equations for Simulating Biomolecular Diffusion-Reaction Processes II: Size Effects on Ionic Distributions and Diffusion-Reaction Rates, Biophys J., {100}, 2011, 2475¨C2485.
%\bibitem{maho00} S. Markus, E. Houstis, A. Catlin, J. Rice, P. Tsompanopoulou, E. Vavalis, D. Gottfried,
%K. Su, and G. Balakrishnan, An agent-based netcentric framework for multidisciplinary
%problem solving environments (MPSE), Internat. J. Comput. Engrg. Sci., {1},
%2000, 33-60.
%\bibitem{maxu95} M. Marion and J. Xu, Error estimates on a new nonlinear Galerkin method based on
%two-grid finite elements, SIAM J. Numer. Anal., {32}, 1995, 1170-1184.
%\bibitem{muxu07} M. Mu and J. Xu, a two-grid method of a mixed Stoke-Darcy model for coupling fluid flow with porous media flow, SIAM J. Numer. Anal., {45}, 2007, 1801-1813.

%\bibitem{quva99} A. Quarteroni and A. Valli, Domain Decomposition Methods for Partial Differential Equations,
%Oxford University Press, Oxford, UK, 1999.


%\%bibitem{sino09} A. Singer and J. Norbury., A Poisson-Nernst-Planck model for biological ion channels - an
%asymptotic analysis in a three-dimensional narrow funnel., SIAM J. Appl. Math., {70}, 2009, 949-968.

%\bibitem{sozh041} Y. H. Song, Y. J. Zhang, C. L. Bajaj, and N. A. Baker, Continuum
%%diffusion reaction rate calculations of wild-type and mutant mouse
%acetylcholinesterase: adaptive finite element analysis, Biophys. J.,
%{3}, 2004, 1558-1566.

%\%bibitem{sozh042} Y. H. Song, Y. J. Zhang, T. Y. Shen, C. L. Bajaj, J. A. McCammon,
%and N. A. Baker, Finite element solution of the steady-state
%Smoluchowski equation for rate constant calculations, Biophys. J.,
%{4}, 2004, 2017-2029.
%\bibitem{susu16} Y. Z. Sun, P. T. Sun, B. Zheng and G. Lin, Error analysis of finite element method for Poisson-Nernst-Planck equations, J. Comput. Appl. Math., {301}, 2016, 28-43.
%\bibitem{whgi11}J. P. Whiteley, K. Gillow, S. J. Tavener and A. C. Walter, Error bounds on block Gauss-Seidel solutions of coupled multiphysics problems, Int. J. Numer. Meth. Engng., {88}, 2011, 1219-1237.
%\bibitem{wusr02} J. Wu, V. Srinivasan, J. Xu, and C. Wang, Newton-Krylov-Multigrid algorithms for battery
%simulation, J. Electrochem. Soc., {149}, 2002, 1342-1348.

%\bibitem{xuzh00}J. Xu and A. Zhou, Local and parallel finite element
%algprithms based on two-grid discretizations,
% Math. Comp., {{69}}, 2000, 881-909.
%\bibitem{xuzh02}J. Xu and A. Zhou, Local and parallel finite element algorithms for eigenvalue problems,
%Acta Mathematicae Applicatea Sinica, English Series, {18}, 2002, 185-200.
\bibitem{a75} R. A. Adams, Sobolev Spaces, Academic Press, New York, 1975.
 \bibitem{bks09} M. Z. Bazant, M. S. Kilic, B. D. Storey and A. Ajdari, Towards an understanding of induced-charge
  electrokinetics
  at large applied voltages in concentrated
solutions, Adv. Colloid Interface Sci., {152}, 2009, 48-88.
 \bibitem{bsdk09} D. S. Bolintineanu, A. Sayyed-Ahmad, H. T. Davis and Y. N. Kaznessis, Poisson-Nernst-Planck models of
  nonequilibrium ion electrodiffusion through a
protegrin transmembrane pore, PLoS Comput. Biol., {5}, 2009, e1000277.
\bibitem{cc65} H. Cohen and J. W. Cooley, The numerical solution of the time-dependent Nernst-Planck equations, Biophys. J., {5},
 1965, 145-162.
 \bibitem{chenchen13} Y. Y. Chen, L. P. Chen and X. C. Zhang, Two-grid method for nonlinear parabolic equations by expanded mixed finite element methods, Numer. Methods Partial Differential Equations, {29}, 2013, 1238-1256.
\bibitem{cck00} A. E. Cardenas, R. D. Coalson and M. G. Kurnikova, Three-dimensional Poisson-Nernst-Planck theoty studies: influence of membrane electrostatics on gramicidin a channel conductancethe, Biophys. J., {79}, 2000, 80-93.
 \bibitem{ck05} R. D. Coalson and M. G. Kurnikova, Poisson-Nernst-Planck theory approach to the calculation of current through
  biological ion channels, IEEE Trans.
Nanobiosci, {4}, 2005, 81-93.

\bibitem{ehl10} B. Eisenberg, Y. K. Hyon and C. Liu, Energy variational analysis of ions in water and channels: Field theory for
  primitive models of complex ionic fluids,
J. Chem. Phys., {133}, 2010, 104104.
%\bibitem{hszhu13}M. Holst, R. Szypowski and Y. R. Zhu, Two-grid methods for semilinear interface problems,
%Number. Methods Partial Differ. Equ., {29}, 2013, 1729-1748.
%\bibitem{ir02} W. Im and B. Roux, Ion permeation and selectivity of OmpF porin a theoretical study based on molecular dynamics,
 %Brownian dynamics, and continuum
 \bibitem{j96} J. Jerome, Analysis of Charge Transport: A Mathematical Theory and Approximation of Semiconductor
 Models, Springer-Verlag, New York, 1996.
% \bibitem{jin15}  J. C. Jin, a two-grid finite-element method for the nonlinear schr\"odinger equations,
%J. Comput. math., {33}, 2015, 146-157.
\bibitem{liudu15} Y. Liu, Y. W. Du, H. Li, J. C. Li and S. He, A two-grid mixed finite element method for a nonlinear fourth-order reaction-diffusion problem with time-fractional derivative, Comput. Math. Appl., {70}, 2015, 2474-2492.

 \bibitem{lhmz} B. Lu, M. J. Holst, J. A. McCammon and Y.C. Zhou, Poisson-Nernst-Planck equations for simulating biomolecular
  diffusion¨Creaction processes I: finite
element solutions, J. Comput. Phys., {229}, 2010, 6979-6994.
\bibitem{lzhb} B. Lu, Y. Zhou, G. A. Huber, S. D. Bond, M. J. Holst and J. A. McCammon, Electrodiffusion: a continuum modeling framework for biomolecular systems with realistic spationtemporal resolution, J. Chem. Phys.,{127}, 2007, 135102.

\bibitem{m86} P. A. Markowich, The Stationary Semiconductor Device Equation, Springer-Verlag, New York, 1986.

\bibitem{mcr12} J. Marcicki, A. T. Conlisk and G. Rizzoni, Comparison of limiting descriptions of the electrical double layer
 using a simplified lithium-ion battery model,
ECS Trans., {41}, 2012, 9-21.
\bibitem{mm09} S. R. Mathur and J. Y. Murthy, A multigrid method for the Poisson-Nernst-Planck equations, Int. J. Heat Mass
 Transfer., {52}, 2009, 4031-4039.
 \bibitem{n91} J. S. Newman, Electrochemical Systems, Prentice Hall, 1991.

 %\bibitem{szs04} Y. H. Song, Y. J. Zhang, T. Y. Shen, C. L. Bajaj, J. A. McCammon and N. A. Baker, Finite element solution of the
 %steady-state Smoluchowski equation for rate
%constant calculations, Biophys. J., {4}, 2004, 2017-2029.

\bibitem{suzhli16} Y. Z. Sun, P. T. Sun, B. Zheng and G. Lin, Error analysis of
finite element method for Poisson-Nernst-Planck equations, J.
Comput. Appl. Math., {301}, 2016, 28-43.
\bibitem{wuwu} L. Wu, M. B. Allen, A two-grid method for mixed finite-element solution of reaction-diffusion equations, Numer. Methods Partial Differential Equations, {15}, 1999, 317-332.
\bibitem{wsx02} J. Wu, V. Srinivasan, J. Xu and C. Wang, Newton-Krylov-multigrid algorithms for battery simulation, J.
Electrochem. Soc., {149}, 2002, 1342-1348.
\bibitem{xu92} J. Xu, A new class of iterative methods for nonselfadjoint or indefinite problems, SIAM J.
Numer. Anal., {29}, 1992, 303-319.
\bibitem{xu94} J. Xu, A novel two-grid method for semilinear elliptic equations, SIAM J. Sci. Comput., {15},
1994, 231-237.
\bibitem{xu96} J. Xu, Two-grid discretization techniques for linear and nonlinear PDEs, SIAM J. Numer.
Anal., {33}, 1996, 1759-1777.
\bibitem{yalu13}Y. Yang and B. Z. Lu, An Error Analysis for the Finite Element Approximation to the
Steady-state Poisson-Nernst-Planck Equations, Adv. Appl. Math.
Mech., {5}, 2013, 113-130.
\bibitem{yaluxie} Y. Yang, B. Z. Lu and Y. Xie, An Decoupling Two-grid Method for the
Steady-state Poisson-Nernst-Planck Equations, submitted.

%\bibitem{zsh15}L. Y. Zuo, P. T. Sun and Y. R. Hou, A Two-grid Method decoupling method for the
%mixed Stokes-Darcy model, J.
%Comput. Appl. Math., {275}, 2015, 139-147.
%
%\bibitem{yazh06}Y. Yang and A. Zhou£¬Two-scale finite element Green¡¯s function approximations with applications
%to electrostatic potential computation, J. Syst. Sci. Complex., {23}, 2010, 177¨C193.

%\bibitem{zhch11}Q. Zheng, D. Chen, and G. W. Wei, Second-order poisson-nernst-planck solver for ion transport,
%J. Comput. Phys, {230}, 2011, 5239-5262.
%\bibitem
%\bibitem{zhlu08} Y. C. Zhou, B. Z. Lu, G. A. Huber, M. J. Holst
%and J. A. McCammon, Continuum simulations of acetylcholine
%consumption by acetylcholinesterase-a Poisson-Nernst-Planck
%approach, J. Phys. Chem. B, {112} (2), 2008, 270-275.

%%%%%%%%%%%%%%%%%%%%%%%%%%%%%%%%%%55555 no use!!!
%\bibitem{ad75}R.A. Adams, Sobolev Spaces, Academic Press, New York, 1975.
%
%\bibitem{brsc}S.C. Brenner and L.R. Scott, The Mathematical Theory of Finite Element Methods,
%Springer-Verlag, New York, 1994.
%
%\bibitem{baho00}N. Baker, M. Holst, and F. Wang, Adaptive
%multilevel finite element solution of the Poisson-Boltzmann equation
%II. refinement at solvent-accessible surfaces in biomolecular
%systems, J. Comput. Chem., {\bf 21} (2000), 1343-1352.
%
%\bibitem{camc90}Z. Cai and S.F. McCormick, On the accuracy of the
%finite volume method for diffusion equations on composite grids,
%SIAM J. Numer. Anal., {\bf 27} (1990), 636-655.
%
%%\bibitem{chho07} L. Chen, M. Holst, and J. Xu, The finite element approximation of
%%the nonlinear Poisson-Boltzmann equation, SIAM J. Numer. Anal., {\bf
%%45} (2007), 2298-2320.
%
%\bibitem{cili91}P.G. Ciarlet and J.L. Lions, Handbook of Numerical
%Analysis, Vol.II, Finite Element Method, Part I , North-Holland,
%Amsterdam, 1991.
%
%\bibitem{cofr972} C.M. Cortis and R.A. Friesner, Numerical solution of
%the Poisson-Boltzmann equation using tetrahedral finite-element
%meshes, J. Comput. Chem., {\bf 18} (1997), 1591-1608.
%
%%\bibitem{dari-duran-padra00} E. Dari, R.G. Duran, and C. Padra,
%%Maximum norm error estimators for three-dimensional elliptic
%%problems, SIAM J. Numer. Anal., {\bf 37} (2000), 683-700.
%
%\bibitem{davis-mccammon90} M.E. Davis and J.A. McCammon,
%Electrostatics in biomolecular structure and dynamics, Chem. Rev.,
%{\bf 90} (1990), 509-521.
%
%\bibitem{ewla91}R.E. Ewing, R.D. Lazarov and P.S. Vassilevski, Locla
%refinement technique for elliptic problem on cell-centered grids, I:
%Error analysis, Math. Comput.,  {\bf 56} (1991) 437-462.
%
%\bibitem{fobr02}F. Fogolari, A. Brigo, and H. Molinari,
%The Poisson-Boltzmann equation for biomolecular electrostatics: a
%tool for structural biology, J. Mol. Recognit., {\bf 15} (2002),
%377-392.
%
%\bibitem{ho94}M. Holst, The Poisson-Boltzmann Equation, Analysis
%and Multilever Numerical Solution, Thesis, UCSD, 1994.
%
%\bibitem{hoba00}M. Holst, N. Baker, and F. Wang, Adaptive
%multilevel finite element solution of the Posson-Boltzmann equation
%I. Algorithms and examples, J. Comput. Chem., {\bf 21} (2000),
%1319-1342.
%
%\bibitem{nitsche79} J. Nitsche, $L^{\infty}-$error analysis for
%finite elements, in: Whiteman, Jr., ed., The Mathematics of Finite
%Elements and Applications, Academic Press, New York, 1979, 173-186.
%
%\bibitem{rasc82}R. Rannacher and R. Scott, Some Optimal Error Estimates
%for Piecewise Linear Finite Element Approximations,  Math. Comput.,
%{\bf 38} (1982), 437-445.
%
%\bibitem{scwa77}A.H. Schatz and L.B. Wahlbin, Interior maximum
%norm estimates for finite element methods, Math. Comput., {\bf 31}
%(1977), 414-442.
%
%\bibitem{sharp91} K.A. Sharp, Incorporating solvent and ion screeing
%into molecular dynamics using the finite-difference Poisson-Boltzman
%method, J. Comput. Chem., {\bf 12} (1991), 454-468.
%
%\bibitem{shho90}K.A. Sharp and B. Honig, Calculating total
%electrostatic energies with the non-linear Poisson-Boltzmann
%equation, J. Phys. Chem., {\bf 94} (1990), 7684-7692.
%
%\bibitem{sharp-honig90}K.A. Sharp and B. Honig, Electrostatic
%interactions in macomolecules: Theory and applications, Annu. Rev.
%Biophys. Biophys. Chem., {\bf 19} (1990), 301-332.
%
%\bibitem{simonson01} T. Simonson, Macomolecular electrostatics:
%Continuum models and their growing pains, Curr. Opin. Struct. Biol.,
%11 (2001), 243-252.
%
%\bibitem{vape92}P.S. Vassilevski, S.I. Petrova, and R.D. Lazarov,
%Finite difference scheme on triangular cel-centered grids with local
%refinement, SIAM J. Sci. Stat. Comput.,  {\bf 13} (1992), 1287-1313.
%
%
%
%\bibitem{yang-zhou08} Y. Yang and A. Zhou, A finite element recovery approach
%to Green's function approximations with applications to
%electrostatic potential computation, J. Comput. Applied Math.,
%{\bf 225}(2009),202-212.
%
%




\end{thebibliography}
\end{document}